\documentclass[a4paper,fleqn]{cas-dc}

\usepackage[authoryear,longnamesfirst]{natbib}

\usepackage[ruled,vlined,linesnumbered]{algorithm2e}
\usepackage{enumerate}
\usepackage{graphicx}
\usepackage{amsmath}
\usepackage{amssymb}
\usepackage{xcolor}
\usepackage{braket}
\usepackage{diagbox}
\usepackage{dsfont}

\ExplSyntaxOn
        	\newcommand{\complexnum}[2]{
    	\str_if_eq:eeTF { \str_head:n {#2} } { - }
     		{#1 - \rmi \, \str_tail:n {#2}}
            	{#1 + \rmi \, #2}
        	}
\ExplSyntaxOff

\usepackage{listings}
\usepackage{lstautogobble}

\newcommand{\rmi}{\mathrm{i}}

\newcommand{\opH}{\mathbb{H}}
\newcommand{\opK}{\mathbb{K}}
\newcommand{\opV}{\mathbb{V}}

\newcommand{\Id}{\mathds{1}}
\newcommand{\dt}{\delta t}
\newcommand{\dx}{\delta x}
\newcommand{\Uq}{\tilde{U}_q}

\newtheorem{theorem}{Theorem}

\newtheorem{proposition}[theorem]{Proposition}

\newdefinition{rmk}{Remark}
\newproof{pf}{Proof}

\begin{document}
\let\WriteBookmarks\relax
\def\floatpagepagefraction{1}
\def\textpagefraction{.001}

\shorttitle{High-order Crank-Nicolson method}
\shortauthors{L.~Freschi et al.}

\title[mode=title]{High-order Crank-Nicolson method via diagonal
  Pad\'{e} factorization}

\author[1]{Lou Freschi}[orcid=0000-0002-2460-295X]
\credit{Conceptualization, Methodology, Software, Formal Analysis, Investigation, Writing -- Review \& Editing}

\author[1]{Camille L\'{e}v\^{e}que}[orcid=0000-0002-7087-1798]
\credit{Methodology, Validation, Writing -- Review \& Editing}

\author[1]{J\'{e}r\'{e}mie Caillat}[orcid=0000-0002-1975-6356]
\credit{Validation, Writing -- Review \& Editing}

\author[1]{Richard Ta\"{i}eb}[orcid=0000-0001-6301-6008]
\credit{Conceptualization, Formal Analysis, Validation, Writing -- Review \& Editing}

\author[1]{Jonathan Dubois}[orcid=0000-0002-6976-5488]
\cormark[1]
\ead{jonathan.dubois@sorbonne-universite.fr}
\credit{Conceptualization, Methodology, Software, Formal Analysis, Supervision, Validation, Writing -- Original Draft}

\cortext[1]{Corresponding author}

\affiliation[1]{%
  organization={Sorbonne Universit\'{e}, CNRS,
    Laboratoire de Chimie Physique -- Mati\`{e}re et Rayonnement, LCPMR},
  postcode={75005},
  city={Paris},
  country={France}%
}

\begin{abstract}
We present a numerical method to solve the time-dependent Schr\"{o}dinger equation at arbitrary order of accuracy in the time step, particularly well suited for tridiagonal Hamilton operators. 
The method is based on a factorization of the diagonal Pad\'{e} approximants of the matrix exponential into products of linear factors, which preserves tridiagonal matrix structure at every time step and yields a propagator of order $2q$ in the time step at a cost of $O(qN)$ operations. 
The propagator is exactly unitary by construction, and at lowest order $q=1$ it reduces to the standard Crank-Nicolson algorithm. 
We further extend the method to any Hamiltonian admitting a decomposition into a sum of tridiagonal or diagonal operators via high-order split-operator techniques~\citep{McLachlan2022}, retaining $O(qN)$ scaling at each time step. 
We also analytically show and numerically verify that the precision is dictated by the split-operator method.
This covers a broad range of physically relevant cases, including higher-dimensional systems, higher-order spatial
discretizations, and non-separable potentials.
\end{abstract}

\maketitle

\section{Introduction}
\label{sec:intro}

Since its original formulation for evaluating numerically heat conduction~\citep{Crank1947}, the Crank-Nicolson method has become a fundamental numerical scheme for solving parabolic partial differential equations. Today, Crank-Nicolson algorithms are found in a wide array of scientific and engineering disciplines due to their unconditional stability and second-order accuracy in both time and space. Beyond classical thermodynamics, they are heavily utilized in quantitative finance to solve the Black-Scholes equation for robust option pricing~\citep{Wilmott1993, Duffy2006}, and in computational fluid dynamics to compute viscous diffusion flows~\citep{Anderson1995}. The method is also a standard tool for modeling complex reaction-diffusion systems in mathematical biology~\citep{Murray2002}, as well as in implicit finite-difference time-domain schemes for accelerating electromagnetic simulations~\citep{Taflove2005}. Here, we will focus on the time dependent Schr\"{o}dinger equation (TDSE) as an example~\citep{Tannor2007}.

\subsection{Time-dependent Schr\"{o}dinger equation}

The TDSE arises across a wide range of applications in modern physics and physical chemistry~\citep{Kulander1991}.
It governs the time evolution of the quantum state $\ket{\psi(t)}$ of a system described by a time-independent Hamiltonian $\opH$. 
Setting $\hbar = 1$, its formal solution over a time interval $[0,T]$ reads
\begin{equation}
	\label{eq:formal_sol_TDSE}
 	 \ket{\psi(T)} =
  	\exp\!\bigl(-\rmi \,\opH\, T \bigr)\,\ket{\psi(0)},
\end{equation}
where $\ket{\psi(0)}$ is the initial state.
For clarity we restrict attention to the time-independent case throughout, although the extension to time-dependent Hamiltonians $\opH(t)$ is straightforward via a Magnus expansion or a piecewise-constant approximation of $\opH(t)$ over each time step~\citep{Blanes2009}.

Once $\opH$ is represented as a matrix, evaluating \eqref{eq:formal_sol_TDSE} requires computing the action of a matrix exponential on a vector. 
A broad and physically important class of Hamiltonians admits a tridiagonal matrix representation, including finite-difference discretizations of the kinetic energy operator on a spatial grid, collective spin models, bosonic and fermionic tight-binding chains, and nearest-neighbour lattice models~\citep{Golub2013}. 
More generally, many systems lead to banded matrices, or to Hamiltonians that decompose naturally into a sum of tridiagonal and diagonal contributions.
The present work places its primary focus on tridiagonal $\opH$, and shows how the method extends naturally to banded matrices and, via operator splitting, to the broader class of Hamiltonians admitting such a decomposition.

For long propagation times the interval $[0,T]$ is split into $N_t$ steps of size $\dt = T / N_t$, and the evolution of the quantum state over one time-step interval reads as
\begin{equation}
	\label{eq:main}
 	 \ket{\psi(t+\dt)} = \exp\!\bigl(-\rmi\,\opH \, \dt \bigr)\,\ket{\psi(t)} ,
\end{equation}
with $t=n \, \dt$ for $n=0,1, \ldots , N_t - 1$.
This partitioning strategy is closely related to the scaling and squaring approach~\citep{Higham2005}, although the two settings differ fundamentally: 
here $N_t$ is determined by the required physical time resolution rather than by numerical optimisation of a single matrix exponential computation.
We denote the exact propagator over one time step by
\begin{equation}
    \label{eq:Udt}
    U(\dt) = \exp\!\bigl(-\rmi\,\dt\,\opH\bigr),
\end{equation}
so that Eq.~\eqref{eq:main} reads $\ket{\psi(t+\dt)} = U(\dt)\ket{\psi(t)}$, or in vector notation
\begin{equation}
    \label{eq:uv}
    \mathbf{u} = U(\dt)\,\mathbf{v} ,
\end{equation}
where $\mathbf{u} = \ket{\psi(t+\dt)}$ and $\mathbf{v} = \ket{\psi(t)}$.
The goal is to compute an approximation $\tilde{U}(\dt)$
to $U(\dt)$ satisfying
\begin{equation}
    \label{eq:approx}
    \bigl\|\tilde{U}(\dt) - U(\dt)\bigr\| = O\!\left(\dt^{r+1}\right),
\end{equation}
where $r$ denotes the order of accuracy, such that $\tilde{U}(\dt)\mathbf{v}$ can be computed efficiently for tridiagonal $\opH$.
Different approaches are employed in the literature to construct such approximations, as reviewed below.

\subsection{Existing approaches}

\paragraph{FFT split-operator methods.}
High-order split-operator schemes decompose the propagator into alternating kinetic and potential factors evaluated via the Fast Fourier Transform~\citep{McLachlan2022}. 
They achieve $O(N\log N)$ cost per step and preserve unitarity for chosen splitting coefficients, but require a uniform periodic grid and a separable Hamiltonian.

\paragraph{Krylov subspace methods.}
The action $U(\dt) \mathbf{v}$ can be approximated via Lanczos or Arnoldi iteration in a Krylov subspace~\citep{Park1986}. 
The \texttt{expmv} algorithm of~\citep{AlMohy2011} provides rigorous error control for general sparse matrices, but does not guarantee exact unitarity and may require many Krylov vectors for tight accuracy requirements.

\paragraph{Crank-Nicolson method.}
The simplest scheme with exact norm conservation is the Crank-Nicolson (CN) method~\citep{Crank1947},
\begin{equation}
	\label{eq:CN}
  	\Bigl(\mathds{1} + \rmi\tfrac{\dt}{2}\opH\Bigr)\mathbf{u} = \Bigl(\mathds{1} - \rmi\tfrac{\dt}{2}\opH\Bigr)\mathbf{v},
\end{equation}
which is second-order accurate in $\dt$ and exactly unitary for Hermitian $\opH$. 
When $\opH$ is tridiagonal, both operations cost $O(N)$ via the Thomas algorithm. 
However, CN is limited to second-order accuracy and no established high-order generalization that retains $O(N)$ tridiagonal scaling exists in the literature.

\subsection{Present work}
The goal of this work is to preserve the $O(N)$ efficiency of tridiagonal matrix operations while extending temporal accuracy to arbitrary order. 
We achieve this through two ideas. 
First, we employ the diagonal $[q,q]$ Pad\'{e} approximant~\citep{Pade1892} $R_{q,q}(z) = P_q(z)/P_q(-z)$, which approximates $e^z$ to order $2q$ and $P_q$ is real, ensuring exact unitarity.
Second, we factorize $P_q(z) = \prod_{k=1}^q(1 + \alpha_k z)$, where $\alpha_k$ are the roots of $P_q$. 
Each factor $[\Id - \alpha_k(\rmi\,\dt\,\opH)]$ is tridiagonal whenever $\opH$ is tridiagonal, reducing the full propagator to $q$ sequential Thomas solves and $q$ matrix--vector products, all $O(N)$, for a total cost of $O(qN)$ per time step. 
At $q = 1$ with $\alpha_1 = 1/2$ the method recovers the standard CN scheme~\eqref{eq:CN}.

We further extend the method beyond the single-operator tridiagonal case by combining it with the high-order split-operator framework of~\citep{McLachlan2022}. 
Any Hamiltonian decomposable as $\opH = \sum_i \opH_i$, with each $\opH_i$ tridiagonal or diagonal, can be handled at the same $O(qN)$ cost per time step, covering multi-dimensional systems, higher-order spatial discretizations, and non-separable potentials alike.

The specific contributions of this paper, detailed in the sections that follow, are: 
Section~\ref{sec:method} introduces the diagonal Pad\'{e} propagator, its factorization, and validates it on a one-dimensional benchmark (Section~\ref{sec:numerics_1D}).
Section~\ref{sec:extension} extends the method to higher-dimensional systems via operator splitting and validates the combined algorithm on a two-dimensional system
with a non-separable potential (Section~\ref{sec:numerics_2D}).
Section~\ref{sec:conclusions} concludes.

\section{The diagonal Pad\'{e} propagator and its factorization}
\label{sec:method}

\subsection{Pad\'{e} approximants of $\exp(z)$ and the $R_{q,q}$ family}

The Pad\'{e} approximant of order $[p,q]$ of the exponential is the rational function $R_{p,q}(z) = P_{p,q}(z)/Q_{p,q}(z)$, where $P_{p,q}$ and $Q_{p,q}$ are polynomials of degree $p$ and $q$ chosen so that $\exp(z) = R_{p,q}(z) + O(z^{p+q+1})$~\citep{Pade1892,Higham2005}. 
In the present work we restrict attention to the \emph{diagonal} family $R_{q,q}$, for which $p = q$. 
A key property of this family is that the denominator is entirely determined by the numerator via $Q_{q,q}(z) = P_{q,q}(-z)$, so that the rational approximant takes the form
\begin{equation}
    \label{eq:Rqq}
    R_{q,q}(z) = \frac{P_q(z)}{P_q(-z)},
\end{equation}
where we adopt the shorthand $P_q \equiv P_{q,q}$, defined by
\begin{equation}
    P_q(z) = \sum_{i=0}^{q}
        \frac{(2q-i)!\,q!}{(2q)!\,(q-i)!\,i!}\,z^{i}.
    \label{eq:Ppq}
\end{equation}
Applied to the propagation problem~\eqref{eq:uv}, this gives $U(\dt)= R_{q,q}(-\rmi\,\dt\,\opH) + O(\dt^{2q+1})$, achieving order $2q$ accuracy from a single polynomial $P_q$ of degree $q$. 
The diagonal Pad\'{e} polynomials $P_q$ are listed in Table~\ref{tab:Pade} for $q = 1, \ldots, 5$.

\begin{table*}
    	\centering
    	\footnotesize
    	\setlength{\extrarowheight}{4pt}
    	\caption{
    		Diagonal Pad\'{e} polynomials $P_q(z) \equiv P_{q,q}(z)$ defined by Eq.~\eqref{eq:Ppq}, their coefficients $\alpha_q^{(k)}$ defined by Eq.~\eqref{eq:alpha_roots} (16 decimal digits), and the resulting approximation order $2q$. 
		Complex coefficients appear in conjugate pairs. 
		At $q=1$ the method reduces to the standard CN scheme.
	}
    	\label{tab:Pade}
    	\begin{tabular}{|c|l|c|l|c|}
        \hline
        $q$ & \multicolumn{1}{c|}{$P_q(z)$} & $k$
            & \multicolumn{1}{c|}{$\alpha_q^{(k)}$}
            & Order \\
        \hline
        $1$
          & $1 + \dfrac{1}{2}z$
          & $1$
          & $\dfrac{1}{2}$
          & $2$ \\
        \hline
        \multirow{2}{*}{$2$}
          & \multirow{2}{*}{$1 + \dfrac{1}{2}z + \dfrac{1}{12}z^2$}
          & $1$
          & $\dfrac{1}{4} \left( 1 + \dfrac{\rmi}{\sqrt{3}} \right)$
          & \multirow{2}{*}{$4$} \\
          &
          & $2$
          & $\dfrac{1}{4} \left( 1 - \dfrac{\rmi}{\sqrt{3}} \right)$
          & \\
        \hline
        \multirow{3}{*}{$3$}
          & \multirow{3}{*}{$1 + \dfrac{1}{2}z
                              + \dfrac{1}{10}z^2
                              + \dfrac{1}{120}z^3$}
          & $1$
          & $0.2153144231161122$
          & \multirow{3}{*}{$6$} \\
          &
          & $2$
          & $\complexnum{0.1423427884419439}{0.1357999257081538}$
          & \\
          &
          & $3$
          & $\complexnum{0.1423427884419439}{-0.1357999257081538}$
          & \\
        \hline
        \multirow{4}{*}{$4$}
          & \multirow{4}{*}{$1 + \dfrac{1}{2}z
                              + \dfrac{3}{28}z^2
                              + \dfrac{1}{84}z^3
                              + \dfrac{1}{1680}z^4$}
          & $1$
          & $\complexnum{0.1584337597342824}{0.0474410125711084}$
          & \multirow{4}{*}{$8$} \\
          &
          & $2$
          & $\complexnum{0.1584337597342824}{-0.0474410125711084}$
          & \\
          &
          & $3$
          & $\complexnum{0.0915662402657176}{0.1156626130131276}$
          & \\
          &
          & $4$
          & $\complexnum{0.0915662402657176}{-0.1156626130131276}$
          & \\
        \hline
        \multirow{5}{*}{$5$}
          & \multirow{5}{*}{$1 + \dfrac{1}{2}z
                              + \dfrac{2}{15}z^2
                              + \dfrac{1}{60}z^3
                              + \dfrac{1}{840}z^4
                              + \dfrac{1}{30240}z^5$}
          & $1$
          & $0.1371088129651923$
          & \multirow{5}{*}{$10$} \\
          &
          & $2$
          & $\complexnum{0.1174272543596986}{0.0610497038185382}$
          & \\
          &
          & $3$
          & $\complexnum{0.1174272543596986}{-0.0610497038185382}$
          & \\
          &
          & $4$
          & $\complexnum{0.0640183391577052}{0.0983410691731092}$
          & \\
          &
          & $5$
          & $\complexnum{0.0640183391577052}{-0.0983410691731092}$
          & \\
        \hline
        $\vdots$ & $\vdots$ & $\vdots$ & $\vdots$ & $\vdots$ \\
        \hline
    \end{tabular}
\end{table*}

\subsection{Factorization and root computation}

The polynomial $P_q(z)$ can be factorized as
\begin{equation}
	\label{eq:Ppq_alphak}
	P_q(z) = \prod_{k=1}^q \Bigl[1 + \alpha_q^{(k)}\,z\Bigr], \qquad q \geq 1,
\end{equation}
where the $\alpha_q^{(k)} \in \mathbb{C}$ are the negated reciprocals of the roots $z_q^{(k)}$ of $P_q$,
\begin{equation}
    	\label{eq:alpha_roots}
    	\alpha_q^{(k)} = -\frac{1}{z_q^{(k)}},
\end{equation}
computed to arbitrary precision by any standard polynomial root-finding algorithm.
The factorized form immediatly shows that $P_q(0)=1$, consistent with Eq.~\eqref{eq:Ppq}. 
The following \textit{Mathematica} code computes the $\alpha_q^{(k)}$ to 16 decimal digits for any order $q$:

\begin{lstlisting}[
  language=Mathematica,
  caption={%
    \textit{Mathematica} code for computing the Pad\'{e} roots
    $\alpha_q^{(k)}$ to 16 decimal digits. The output matches
    Table~\ref{tab:Pade}.},
  label=lst:mathematica]
NOrder    = 5;   (* Pade order q *)
NAccuracy = 16;  (* number of decimal digits *)

coeffs = Table[
  Binomial[NOrder, j] / (Binomial[2*NOrder, j] * j!),
  {j, 1, NOrder}
];

poly   = 1 + Sum[coeffs[[j]] * z^j, {j, 1, NOrder}];
roots  = z /. NSolve[poly == 0, z, NAccuracy];
values = N[-1/roots, NAccuracy];

values // MatrixForm
\end{lstlisting}

The $\alpha_q^{(k)}$ satisfy two general properties that are useful both for verification and for the unitarity proof below. 
First, since $P_q$ has real coefficients, complex roots appear in conjugate pairs,
\begin{equation}
    \label{eq:conj_pair}
    \alpha_q^{(k+1)} = \bigl[\alpha_q^{(k)}\bigr]^{\star},
\end{equation}
so that for odd $q$ there is exactly one real root and $(q-1)/2$ conjugate pairs, while for even $q$ all $q/2$ roots are complex conjugate pairs. 
Second, identifying the coefficient of $z^1$ in Eqs.~\eqref{eq:Ppq} and~\eqref{eq:Ppq_alphak} gives
\begin{equation}
    \label{eq:alpha_sum}
    \sum_{k=1}^q \alpha_q^{(k)} = \frac{1}{2},
\end{equation}
which provides a useful numerical check: For $q=1$, $\alpha_1^{(1)} = 1/2$ exactly, while for $q=2$, the two roots $\alpha_2^{(1,2)} = (1 \pm \rmi/\sqrt{3})/4$ obviously satisfy Eq.~\eqref{eq:alpha_sum} exactly. 
The values of $\alpha_q^{(k)}$ for $q = 1,\ldots,5$ are listed in Table~\ref{tab:Pade}.

\subsection{Factored Padé propagator \label{sec:Pade_propagator}}

With the factorization~\eqref{eq:Ppq_alphak} in hand, the propagation step $\mathbf{u} \approx U (\dt) \, \mathbf{v}$ can be written in high-order CN form
\begin{equation}
    \label{eq:CN_high_order}
    \prod_{k=1}^q \Bigl[\Id + \alpha_q^{(k)} ( \rmi\,\dt\,\opH ) \Bigr]
    \,\mathbf{u}
    =
    \prod_{k=1}^q \Bigl[\Id - \alpha_q^{(k)} ( \rmi\,\dt\,\opH ) \Bigr]
    \,\mathbf{v},
\end{equation}
or equivalently $P_q(\rmi\,\dt\,\opH)\,\mathbf{u} = P_q(-\rmi\,\dt\,\opH)\,\mathbf{v}$.
For $q=1$, this reduces exactly to the standard CN equation~\eqref{eq:CN}.
For higher values of $q$, the right-hand side is evaluated first as a sequence of $q$ tridiagonal matrix--vector multiplications, then the left-hand side is solved as a sequence of $q$ tridiagonal linear systems via the Thomas algorithm.
Thus, it remains identical in structure to CN but extended to order $2q$ at a cost of $O(qN)$ operations. 
The formal properties of this propagator are established in Theorem~\ref{thm:main} below, using the compact symmetric form $\Uq = [P_q(\rmi\,\dt\,\opH)]^{-1}P_q(-\rmi\,\dt\,\opH)$ which follows from the commutativity of all factors.

\begin{theorem}
\label{thm:main}
Let $\opH \in \mathbb{C}^{N\times N}$ be a Hermitian tridiagonal
matrix, let $\dt \in \mathbb{R}$ (the time step may be negative), and let $\alpha_q^{(1)}, \ldots,
\alpha_q^{(q)}$ be defined by Eq.~\eqref{eq:alpha_roots} from
the roots of $P_q$ as given by Eq.~\eqref{eq:Ppq_alphak}.
Since all factors $\bigl[\Id \pm \alpha_q^{(k)}\rmi\,\dt\,\opH
\bigr]$ are functions of the same operator $\opH$, they mutually
commute, and the propagator~\eqref{eq:CN_high_order} can be
written in the symmetric form
\begin{equation}
    \label{eq:Uq}
    \Uq (\dt)   =
    \Bigl[P_q\bigl(\rmi\,\dt\,\opH\bigr)\Bigr]^{-1}
    P_q\bigl(-\rmi\,\dt\,\opH\bigr),
\end{equation}
where $P_q$ is defined by Eq.~\eqref{eq:Ppq_alphak}. This
propagator satisfies the following properties:
\begin{enumerate}[(i)]
    \item \textbf{Exact unitarity:} $\Uq^\dagger \Uq = \Id$;
    \item \textbf{Accuracy:}
          $\bigl\| \Uq - U \bigr\|
          = O\!\left(\dt^{2q+1}\right)$;
    \item \textbf{Linear complexity:} 
    	Each factor $\bigl[\Id \pm \alpha_q^{(k)}\rmi\,\dt\,\opH\bigr]$ is tridiagonal, so that $\Uq\,\mathbf{v}$ is obtained from $q$ tridiagonal matrix--vector products and $q$ tridiagonal solves (Thomas algorithm), each costing $O(N)$ operations, i.e.\ $O(qN)$ operations in total.
\end{enumerate}
\end{theorem}

\begin{pf}
We prove each property in turn.

\textbf{(i) Exact unitarity.}
Since $P_q$ has real coefficients and $\opH$ is Hermitian,
\begin{equation*}
    \Bigl[P_q\bigl(-\rmi\,\dt\,\opH\bigr)\Bigr]^\dagger
    = P_q\bigl(\rmi\,\dt\,\opH\bigr),
\end{equation*}
and therefore
\begin{eqnarray*}
    \Uq^\dagger &=&
    \Bigl[P_q\bigl(-\rmi\,\dt\,\opH\bigr)\Bigr]^\dagger
    \Bigl[P_q\bigl(\rmi\,\dt\,\opH\bigr)\Bigr]^{-\dagger} \\
    &=& P_q\bigl(\rmi\,\dt\,\opH\bigr) 
      \Bigl[P_q\bigl(-\rmi\,\dt\,\opH\bigr)\Bigr]^{-1} = \Uq^{-1},
\end{eqnarray*}
so $\Uq^\dagger \Uq = \Id$. \qed

\textbf{(ii) Accuracy.}
Setting $z = -\rmi\,\dt\,\opH$, the form~\eqref{eq:Uq} gives $\Uq = P_q(-z)^{-1}P_q(z) = R_{q,q}(z) = R_{q,q}(-\rmi\,\dt\,\opH)$ (cf.\ Eq.~\eqref{eq:Rqq}).
By definition of the diagonal Pad\'{e} approximant, $\exp(z) = R_{q,q}(z) + O(z^{2q+1})$~\citep{Higham2005}, and therefore $\bigl\| \Uq -  U \bigr\| = O(\dt^{2q+1})$. \qed

\textbf{(iii) Linear complexity.}
Since $\opH$ is tridiagonal, this structure is conserved in each factor $\Id \pm \alpha_q^{(k)}\rmi\,\dt\,\opH$.
A tridiagonal matrix--vector product costs $O(N)$ operations, and so does a tridiagonal linear solve via the Thomas algorithm. 
The product~\eqref{eq:CN_high_order} consists of $q$ such products and $q$ such solves, so the total cost of $\Uq\,\mathbf{v}$ is $O(qN)$. \qed
\end{pf}

\begin{rmk}
The roots $z_q^{(k)}$ of $P_q(z)$ all have strictly negative real part for any $q \geq 1$, so that $\text{Re}(\alpha_q^{(k)}) > 0$ for all $k$ and $q$ (cf.\ Table~\ref{tab:Pade}). 
This follows from the A-stability of the diagonal Pad\'{e} approximants $R_{q,q}$~\citep{Hairer1996}:
since $|R_{q,q}(z)| \leq 1$ for $\text{Re}\,z \leq 0$, $R_{q,q}$ has no poles in the closed left half-plane, i.e.\ $P_q(-z)$ has no roots there, and therefore all roots of $P_q(z)$ lie in the open left half-plane. 
As a consequence, for any Hermitian $\opH$ with real eigenvalues $\lambda$ and any $\dt > 0$, the eigenvalues $1 + \rmi\,\alpha_q^{(k)}\dt\,\lambda$ of $\Id + \alpha_q^{(k)}\rmi\,\dt\,\opH$ satisfy $|1 + \rmi\,\alpha_q^{(k)}\dt\,\lambda| \geq \text{Re}(\alpha_q^{(k)})/|\alpha_q^{(k)}| > 0$, uniformly in $\dt$ and $\lambda$. 
The same bound holds for every leading principal submatrix, which is of the same form with a Hermitian block of $\opH$.
Hence every pivot of the Thomas algorithm is non-zero and the algorithm never breaks down, for any real $\dt$. 
Together with the exact unitarity of Theorem~\ref{thm:main}(i), this makes the method unconditionally stable with respect to $\dt$, exactly as the standard CN scheme~\eqref{eq:CN} to which it reduces at $q=1$.
\end{rmk}

\subsection{Numerical experiments: One-dimensional free particle}
\label{sec:numerics_1D}

We validate the factored Pad\'{e} propagator of Section~\ref{sec:Pade_propagator} on a one-dimensional free particle. 
The reference solution $U(T) \mathbf{v}$ is computed directly at the final time $T$ with the \texttt{expmv} algorithm of~\citet{AlMohy2011}, which evaluates the action of the matrix exponential through a truncated Taylor series with automatic selection of the degree and of the number of sub-steps. 
Since the same discretized Hamiltonian is used for the reference and for the propagator, the comparison isolates the temporal error of our method, down to the accuracy of the reference itself, which is a few $10^{-12}$ here (see below).

We consider a free particle on a uniform grid of $N$ points with spacing $\dx$, whose Hamiltonian is discretized by the standard second-order finite difference Laplacian,
\begin{equation}
    \label{eq:H_1D}
    \opH = -\frac{1}{2\dx^2}
    \begin{pmatrix}
        -2 &  1 &        &   \\
         1 & -2 &  1     &   \\
           & \ddots & \ddots & \ddots \\
           &        &  1 & -2
    \end{pmatrix},
\end{equation}
which is Hermitian and tridiagonal, the ideal setting for the method of Section~\ref{sec:method}. 
The initial state is a Gaussian wavepacket centered at $-5$~a.u. of standard deviation 2 a.u. with initial wave vector $k_x = 2$~a.u. and the total propagation time is $T = 4$~a.u.
The space grid consists of $512$ points centered at the origin with $\dx=0.1$~a.u.

\begin{figure}[t]
	\centering
	\includegraphics[width=.33\textwidth]{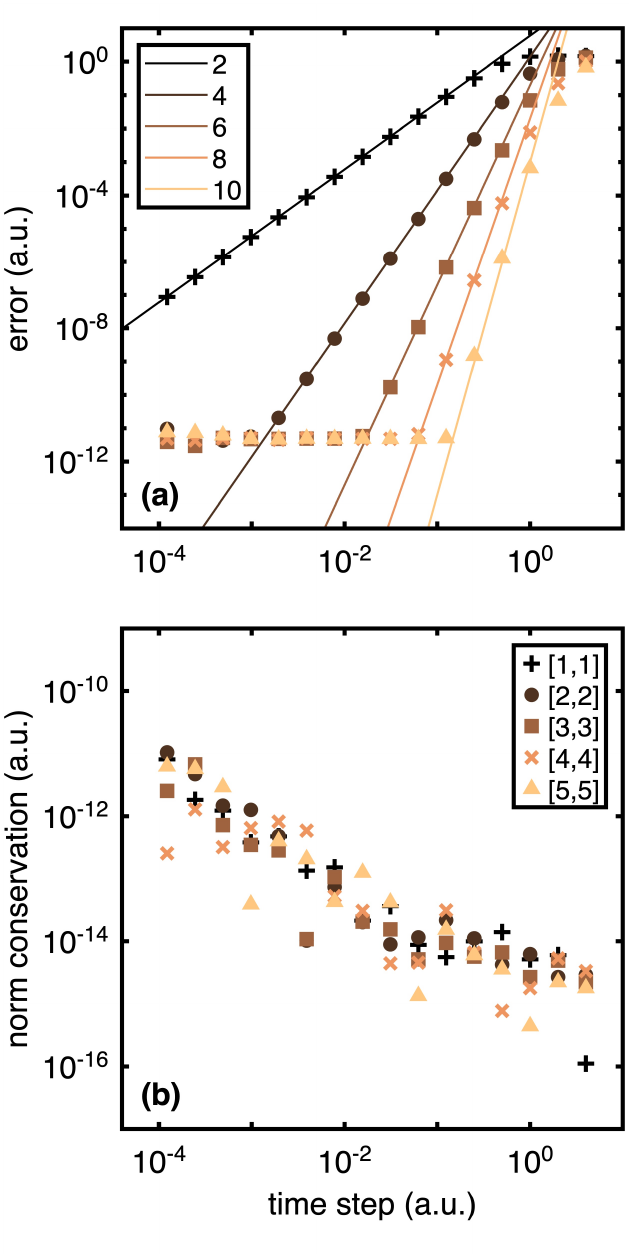}
	\caption{One-dimensional free particle~\eqref{eq:H_1D}.
        (a) $L^2$ error at time $T = 4$~a.u.\ as a function of time step $\dt$ for the diagonal Pad\'{e} families $[q,q]$, $q = 1,\ldots,5$. 
        Solid lines show the theoretical slopes $O(\dt^{2q})$. 
        The error saturates at $\sim 10^{-14}$ (double-precision machine epsilon) and increases slightly for very small $\dt$ due to floating-point accumulation.
        (b) Norm conservation at time $T = 4$~a.u.\ as a function of time step $\dt$.
        All schemes maintain norm deviation at the level of machine precision ($\sim 10^{-14}$) for all $\dt$, confirming exact unitarity by construction (Theorem~\ref{thm:main}(i)).  
        The slight increase for small $\dt$ reflects floating-point accumulation over a large number of steps.}
	\label{fig:error_1D}
\end{figure}

\begin{figure}
	\centering
	\includegraphics[width=.33\textwidth]{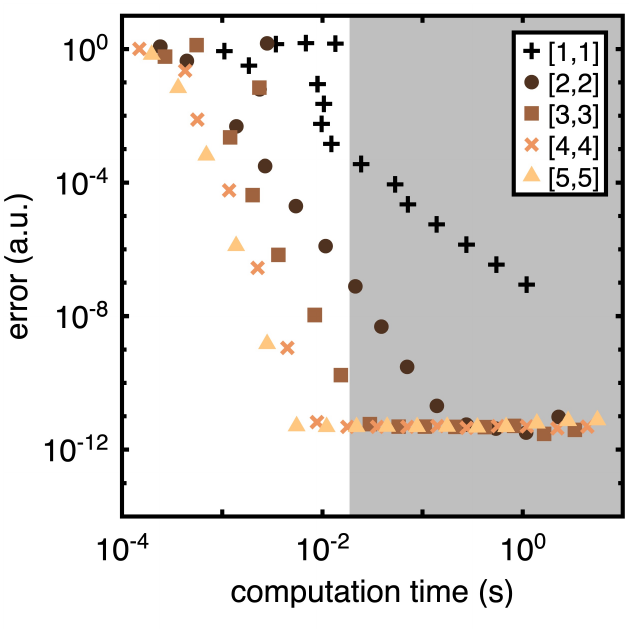}
	\caption{$L^2$ error at time $T = 4$~a.u.\ as a function of the wall-time computation for the diagonal Pad\'{e} families $[q,q]$, $q = 1,\ldots,5$. 
	The gray area indicates the time above which the computation takes longer than \texttt{expmv} function (wall-time computation: 0.018648 s).}
	\label{fig:err_vs_time}
\end{figure}

We apply the factored Pad\'{e} propagator for $q = 1, 2, 3, 4, 5$, corresponding to the diagonal families $[q,q]$ with approximation orders $2, 4, 6, 8, 10$ respectively, over a range of time steps $\dt \in [10^{-4}, 4]$~a.u. 
The $L^2$ error at time $T$ with respect to the \texttt{expmv} reference is shown as a function of $\dt$ in Fig.~\ref{fig:error_1D}(a), and as a function of wall-clock computation time in Fig.~\ref{fig:err_vs_time}. 
The norm conservation $|1 - \| \Uq^{N_t} \mathbf{v} \| |$ at time $T$ is shown in Fig.~\ref{fig:error_1D}(b).

\paragraph{Convergence order.}
Figure~\ref{fig:error_1D}(a) confirms that the error scales as $O(\dt^{2q})$ for each $q$, in agreement with Theorem~\ref{thm:main}(ii). 
The theoretical slopes $2,4,6,8,10$ are shown as solid lines and match the numerical data over several decades. 
For small $\dt$, all curves saturate at a common floor of $\approx 5\times10^{-12}$, independent of both $q$ and $\dt$, and hence of the number of time steps $N_t=T/\dt$, which varies here from $\sim 40$ to $4\times10^{4}$. 
This floor is therefore not a property of the propagator but of the reference solution: \texttt{expmv} evaluates the action of the exponential through a truncated Taylor series requiring a large number of matrix--vector products for the present spectral range ($\|\opH\|\,T \approx 800$), and the accumulated rounding errors limit its accuracy to a few $10^{-12}$~\citep{AlMohy2011}. 
We verified this by comparing with the exact solution obtained from the analytic eigenvectors of the discrete Laplacian~\eqref{eq:H_1D}: 
the \texttt{expmv} reference deviates from it by $5\times10^{-12}$, whereas the $[5,5]$ propagator with $\dt = 0.1$~a.u.\ deviates by only $2\times10^{-13}$.

\paragraph{Norm conservation.}
Figure~\ref{fig:error_1D}(b) shows the norm deviation $|1-\|\Uq^{N_t}\mathbf{v}\||$. 
It increases from $\sim10^{-15}$ for the largest time steps to $\sim10^{-11}$ at $\dt=10^{-4}$~a.u., roughly in proportion to $N_t$, for all $q$. 
Since $\Uq$ is exactly unitary (Theorem~\ref{thm:main}(i)), this deviation is entirely due to floating-point arithmetic: 
each time step introduces a relative error of the order of the machine epsilon $\epsilon\approx1.1\times10^{-16}$, and these errors accumulate linearly over the $N_t$ steps. 
The same accumulation contributes to the error of Fig.~\ref{fig:error_1D}(a), but remains below the floor set by the reference solution, except at the smallest time steps where both reach $\sim10^{-11}$. 
Rounding accumulation thus favours large time steps, and hence high orders: extrapolating the $O(\dt^2)$ slope, the standard CN scheme ($q=1$) would require $N_t\sim10^7$ steps to reach the reference floor, at which point the accumulated rounding error ($\sim N_t\epsilon$) would already exceed $10^{-10}$.

\paragraph{Optimal operating point.}
Combining these observations, the optimal strategy is to use the largest available $q$ with the largest $\dt$ that keeps the error at the desired level. 
For the present experiment, the $[5,5]$ propagator with $\dt\approx0.1$~a.u.\ reaches the accuracy floor of the reference solution with only $N_t=40$ time steps, which also minimizes the accumulated rounding error and the computation time.

\section{Extension to higher-dimensional systems and split-operator
         methods}
\label{sec:extension}


The method of Section~\ref{sec:method} applies directly to any Hermitian tridiagonal matrix $\opH$. 
Many Hamiltonians of interest are not tridiagonal, but can be written as a sum
\begin{equation}
  	\label{eq:H_decomp}
  	\opH = \sum_{i=1}^{S} \opH_i ,
\end{equation}
of Hermitian contributions $\opH_i \in \mathbb{C}^{N \times N}$ that are each tridiagonal, up to a fixed ordering of the grid points, or diagonal. 
On a $d$-dimensional grid, the kinetic energy splits into one operator per direction, each consisting of independent one-dimensional tridiagonal blocks, while the potential is diagonal and its exponential is applied exactly.
Similarly, a $(2M+1)$-point finite-difference Laplacian in one dimension is a sum of $M$ operators coupling sites at distance $m=1,\ldots,M$. 
Since the $\opH_i$ generally do not commute, we combine the factored Pad\'{e} propagator with high-order operator splitting, retaining exact unitarity and $O(qN)$ scaling per time step.

\subsection{The split-operator method}

The naive factorization $U (\dt) = \prod_i \exp(-\rmi\,\dt\,\opH_i)$ fails: the Baker--Campbell--Hausdorff formula gives correction terms proportional to commutators $[\opH_i, \opH_j]$ and higher. 
Following~\citep{McLachlan2022}, we define the elementary maps
\begin{align}
    \chi(h)      &= e^{-\rmi\opH_1 h}\,e^{-\rmi\opH_2 h} \cdots e^{-\rmi\opH_S h}, \label{eq:chi} \\
    \chi^\star(h) &= e^{-\rmi\opH_S h}\cdots e^{-\rmi\opH_2 h} \,e^{-\rmi\opH_1 h}, \label{eq:chi_star}
\end{align}
where $\chi^\star(h) = [\chi(-h)]^{-1}$ denotes the \emph{adjoint method} of $\chi(h)$~\citep{McLachlan1995, Hairer2006, McLachlan2022}.
For a product of operator exponentials, it amounts to reversing their order.
The propagator over one time step is then approximated as
\begin{equation}
	U(\dt) = U_s (\dt) + O (\dt^{r+1}) ,
\end{equation}
by the composition
\begin{equation}
    \label{eq:splitting}
    U_s (\dt) =
    \chi(\beta_1\dt)\,\chi^\star(\beta_2\dt)\,
    \chi(\beta_3\dt)\,\chi^\star(\beta_4\dt)\cdots
    \chi^\star(\beta_{2N_s}\dt),
\end{equation}
where the coefficients $\{\beta_s\}_{s=1}^{N_s}$ are chosen to satisfy the order conditions up to a given order $r$, and the full sequence of $2N_s$ stages is obtained by palindromic reflection
\begin{equation}
    \label{eq:palindrome}
    \beta_{2N_s+1-s} = \beta_s, \qquad s = 1,\ldots,N_s,
\end{equation}
which guarantees time-reversibility and exact unitarity of the full propagator. 
Only the $N_s$ independent coefficients need to be stored.
Table~\ref{tab:splitting} lists them for the Verlet (also known Strang splitting or leapfrog), BM4 and BM6 schemes.

A key advantage of this formulation~\citep{McLachlan2022} is that the coefficients $\{\beta_s\}$ are \emph{independent of the number of operators} $S$ in the decomposition~\eqref{eq:H_decomp}. 
The same table applies whether $S = 2$ (standard two-dimensional problem), $S = 3$ (three-dimensional problem), or any larger $S$, with no re-optimization required. 

\begin{table}
  \centering
  \footnotesize
  \setlength{\extrarowheight}{3pt}
  \caption{Independent splitting coefficients $\beta_s$ for
           $s = 1,\ldots,N_s$, following the $\chi/\chi^\star$
           convention of~\citep{McLachlan2022}. The full sequence
           of $2N_s$ stages is recovered by palindromic
           reflection~\eqref{eq:palindrome}.}
  \label{tab:splitting}
  \begin{tabular}{|c|c|c|c|}
    \hline
    $s$ &
    \multicolumn{1}{c|}{Verlet} &
    \multicolumn{1}{c|}{BM4} &
    \multicolumn{1}{c|}{BM6} \\
    \hline
    $1$  & $\phantom{-}1/2$
         & $\phantom{-}0.0792036964311957$
         & $\phantom{-}0.050262764400392$ \\
    $2$  & \multicolumn{1}{c|}{---}
         & $\phantom{-}0.1303114101821663$
         & $\phantom{-}0.098553683500650$ \\
    $3$  & \multicolumn{1}{c|}{---}
         & $\phantom{-}0.2228614958676077$
         & $\phantom{-}0.314960616927694$ \\
    $4$  & \multicolumn{1}{c|}{---}
         & $-0.3667132690474257$
         & $-0.447346482695478$ \\
    $5$  & \multicolumn{1}{c|}{---}
         & $\phantom{-}0.3246481886897062$
         & $\phantom{-}0.492426372489876$ \\
    $6$  & \multicolumn{1}{c|}{---}
         & $\phantom{-}0.1096884778767498$
         & $-0.425118767797691$ \\
    $7$  & \multicolumn{1}{c|}{---}
         & \multicolumn{1}{c|}{---}
         & $\phantom{-}0.237063913978122$ \\
    $8$  & \multicolumn{1}{c|}{---}
         & \multicolumn{1}{c|}{---}
         & $\phantom{-}0.195602488600053$ \\
    $9$  & \multicolumn{1}{c|}{---}
         & \multicolumn{1}{c|}{---}
         & $\phantom{-}0.346358189850727$ \\
    $10$ & \multicolumn{1}{c|}{---}
         & \multicolumn{1}{c|}{---}
         & $-0.362762779254345$ \\
    \hline
    $N_s$   & $1$ & $6$  & $10$ \\
    Order   & $2$ & $4$  & $6$  \\
    \hline
  \end{tabular}
\end{table}

\subsection{The combined algorithm}

\begin{algorithm}[t]
\DontPrintSemicolon
\SetKwInput{Input}{Input}
\SetKwInput{Output}{Output}
\TitleOfAlgo{Strided Thomas Algorithm for $M$-shifted banded systems}
\Input{
    Main diagonal $b \in \mathbb{C}^N$,\\
    Sub-diagonal $a \in \mathbb{C}^{N-M}$ (entries at positions $M+1, M+2, \ldots, N$),\\
    Super-diagonal $c \in \mathbb{C}^{N-M}$ (entries at positions $1, 2, \ldots, N-M$),\\
    Right-hand side $d \in \mathbb{C}^N$
}
\Output{Solution vector $x \in \mathbb{C}^N$, overwriting $d$}
\BlankLine
\tcp{Allocate modified super-diagonal (stride-aware)}
$c' \leftarrow \mathbf{0} \in \mathbb{C}^{N-M}$\;
\BlankLine
\tcp{Loop over the $M$ independent subsystems, NO re-indexing needed}
\For{$k = 1$ \KwTo $M$}{
    \BlankLine
    \tcp{--- Forward sweep for subsystem $k$, stride $M$ ---}
    $c'_k \leftarrow c_k / b_k$\;
    $d_k \leftarrow d_k / b_k$\;
    \BlankLine
    \For{$i = k+M,\ k+2M,\ \ldots \leq N$}{
        $w        \leftarrow a_{i-M} \cdot c'_{i-M}$\;
        $\text{denom} \leftarrow b_i - a_{i-M} \cdot c'_{i-M}$\;
        \uIf{$i \leq N - M$}{
            $c'_i \leftarrow c_i / \text{denom}$\;
        }
        $d_i \leftarrow (d_i - a_{i-M} \cdot d_{i-M}) / \text{denom}$\;
    }
    \BlankLine
    \tcp{--- Backward substitution for subsystem $k$, stride $M$ ---}
    \tcp{Last index of subsystem $k$}
    $\ell \leftarrow k + M \cdot \lfloor (N - k) / M \rfloor$\;
    $x_\ell \leftarrow d_\ell$\;
    \BlankLine
    \For{$i = \ell - M,\ \ell - 2M,\ \ldots \geq k$}{
        $x_i \leftarrow d_i - c'_i \cdot x_{i+M}$\;
    }
}
\caption{Solves the $M$-shifted banded system $Ax = d$ in $O(N)$
         operations. At $M=1$ this reduces to the standard Thomas
         algorithm.}
\label{alg:strided_thomas}
\end{algorithm}

Each elementary factor $\exp(-\rmi\,\beta_s\,\dt\,\opH_i)$ in Eqs.~\eqref{eq:chi} and~\eqref{eq:chi_star} is the exponential of a scaled Hermitian matrix that is tridiagonal, up to a fixed ordering of the grid points, or diagonal. 
Diagonal contributions, such as the potential, are exponentiated exactly. 
For the tridiagonal ones, we use the factored Pad\'{e} propagator of Section~\ref{sec:method} with effective time step $h=\beta_s\,\dt$, which we denote
\begin{equation}
    	\label{eq:combined}
    	\Uq(h,\opH_i) = \prod_{k=1}^q \Bigl[\Id + \alpha_q^{(k)} ( \rmi\,h\,\opH_i ) \Bigr]^{-1} \Bigl[\Id - \alpha_q^{(k)} ( \rmi\,h\,\opH_i ) \Bigr] ,
\end{equation}
with $\Uq(h,\opH_i) \approx \exp\!\bigl(-\rmi\,h\,\opH_i\bigr)$ so that $\Uq(\dt)$ of Section~\ref{sec:method} corresponds to $\Uq(\dt,\opH)$.
By Theorem~\ref{thm:main}(iii), each elementary propagator $\Uq(h,\opH_i)$ costs $O(qN)$ operations.
This remains true when $\opH_i$ is tridiagonal only up to a permutation of the grid points (Remark~\ref{rmk:strided}).

\begin{rmk}\label{rmk:strided}
In general, an operator $\opH_i$ in Eq.~\eqref{eq:H_decomp} is tridiagonal only up to a fixed permutation of the grid points. 
A frequent case is that of an $M$-shifted tridiagonal matrix, whose non-zero entries lie on the main diagonal and on the diagonals at offsets $\pm M$. 
Such a matrix decouples into $M$ independent tridiagonal systems, formed by the entries $\{k, k+M, k+2M,\ldots\}$ for $k=1,\ldots,M$.
Algorithm~\ref{alg:strided_thomas} solves them directly on the original vector by running the Thomas algorithm on every $M$-th entry (a \emph{strided} access), without copying or reordering the data. 
Its cost is $O(N)$ regardless of $M$, and it reduces to the standard Thomas algorithm at $M=1$. 
The same strided access applies to the tridiagonal matrix--vector products, so that each elementary propagator $\Uq(h,\opH_i)$ costs $O(qN)$ operations.
\end{rmk}

\begin{rmk}
The Thomas algorithm is inherently sequential, so that no parallelism is available for a single tridiagonal system. 
However, whenever an operator $\opH_i$ decouples into $n$ independent tridiagonal systems, as for an $M$-shifted matrix ($n=M$, Remark~\ref{rmk:strided}) or a block-diagonal matrix with tridiagonal blocks, the elementary propagator $\Uq(h,\opH_i)$ splits into $n$ independent problems, each requiring $q$ tridiagonal products and $q$ Thomas solves. 
These problems are embarrassingly parallel and map directly onto shared-memory (OpenMP), distributed-memory (MPI) or GPU (batched Thomas) architectures, while diagonal contributions reduce to elementwise multiplications. 
The sequential work is confined to each individual system, of size $N/n$, and is a negligible fraction of the $O(qN)$ total when $n$ is large. 
This situation is generic on multi-dimensional grids, where each directional operator decouples into one independent system per grid line, and provides an additional incentive to combine the factored Pad\'{e} propagator with operator splitting.
\end{rmk}

Splitting schemes of order $r>2$ contain negative coefficients $\beta_s$ (see Table~\ref{tab:splitting}), so that some elementary propagators are applied with a negative time step $h<0$. 
For the Schr\"{o}dinger equation this is harmless: 
Theorem~\ref{thm:main} holds for any real time step, so each factor $\Uq(\beta_s\dt,\opH_i)$ remains exactly unitary and well defined, whatever the sign of $\beta_s$. 
Negative coefficients would instead be problematic for dissipative problems, such as imaginary-time propagation or Hamiltonians with complex absorbing potentials, for which backward sub-steps amplify the solution.

The full algorithm is given in Algorithm~\ref{alg:combined}. 
Its properties follow directly from Theorem~\ref{thm:main} and the splitting approximation~\eqref{eq:splitting}:

\begin{proposition}\label{prop:combined}
Let $\Uq(h,\opH_i)$ be given by Eq.~\eqref{eq:combined}, replaced by
the exact exponential when $\opH_i$ is diagonal, and define
\begin{align}
  	\tilde{\chi}_q(h)       &= \Uq(h,\opH_1)\,\Uq(h,\opH_2)\cdots\Uq(h,\opH_S),
  	\label{eq:chi_q}\\
  	\tilde{\chi}_q^\star(h) &= \Uq(h,\opH_S)\cdots\Uq(h,\opH_2)\,\Uq(h,\opH_1).
  	\label{eq:chi_q_star}
\end{align}
If the coefficients $\{\beta_s\}$ make the composition~\eqref{eq:splitting} of order $r$, the combined propagator
\begin{equation}
  	\mathcal{U}_q(\dt) = \tilde{\chi}_q(\beta_1\dt)\, \tilde{\chi}_q^\star(\beta_2\dt)\cdots \tilde{\chi}_q^\star(\beta_{2N_s}\dt)
  	\label{eq:combined_prop}
\end{equation}
satisfies, for every $q\geq1$:
\begin{enumerate}[(i)]
  	\item \textbf{Accuracy:} $\mathcal{U}_q(\dt) = e^{-\rmi\dt\opH} + O(\dt^{r+1})$, and the global error after $N_t=T/\dt$ steps is $O(\dt^{r})$;
	\item \textbf{Exact unitarity:} $\mathcal{U}_q^\dagger\,\mathcal{U}_q = \Id$;
  	\item \textbf{Complexity:} the cost per time step is $O(q\,S\,N_s\,N)$, i.e.\ $O(N)$ for fixed $q$, $S$ and $N_s$.
\end{enumerate}
\end{proposition}

\begin{pf}
(i) Since all factors in Eq.~\eqref{eq:combined} commute,
\begin{equation}
  	\Uq(-h,\opH_i) = \bigl[P_q(-\rmi h\opH_i)\bigr]^{-1} P_q(\rmi h\opH_i) = \bigl[\Uq(h,\opH_i)\bigr]^{-1},
  	\label{eq:Uq_symmetric}
\end{equation}
i.e.\ each elementary propagator is time-symmetric, as is the exact exponential. 
Consequently $[\tilde{\chi}_q(-h)]^{-1} = \tilde{\chi}_q^\star(h)$, so that $\tilde{\chi}_q^\star$ is exactly the adjoint of $\tilde{\chi}_q$. 
Moreover, $\tilde{\chi}_q(h) = e^{-\rmi h\opH} + O(h^2)$, so $\tilde{\chi}_q$ is a consistent first-order method. 
Equation~\eqref{eq:combined_prop} is therefore the composition of a first-order method and its adjoint, with unchanged coefficients $\{\beta_s\}$. 
The order conditions of such compositions involve only the $\beta_s$ and are independent of the underlying method~\citep{McLachlan1995,Hairer2006}. 
This is precisely what makes the coefficients of Table~\ref{tab:splitting} valid for any $S$~\citep{McLachlan2022}. 
The order $r$ is thus preserved. 
Since $\mathcal{U}_q$ is unitary by (ii), local errors accumulate at most linearly over the $N_t$ steps, giving a global error $O(\dt^r)$. \qed

(ii) Each factor $\Uq(\beta_s\dt,\opH_i)$ is exactly unitary for any real time step by Theorem~\ref{thm:main}(i), and so is the exact exponential of a diagonal Hermitian matrix; a product of unitary matrices is unitary. \qed

(iii) A time step consists of $2N_s$ stages, each applying $S$ elementary propagators of cost $O(qN)$ (Remark~\ref{rmk:strided}). \qed
\end{pf}

\begin{rmk}
The Pad\'{e} order $q$ thus enters only through the error constant.
Since $R_{q,q}(-z) = 1/R_{q,q}(z)$, the function $\log R_{q,q}$ is odd,
\begin{equation}
  \log R_{q,q}(z) = z + (-1)^{q+1}\frac{(q!)^2}{(2q)!\,(2q+1)!}\,z^{2q+1}
                    + O(z^{2q+3}),
  \label{eq:logPade}
\end{equation}
and, for commuting $\opH_i$, the Pad\'{e} deviations add up in $\log\mathcal{U}_q$ with weight $\sum_{s=1}^{2N_s}\beta_s^{2q+1}$. 
This sum vanishes whenever $2q+1<r$, since $\sum_s\beta_s^{3}=0$ (for $r\ge4$) and $\sum_s\beta_s^{5}=0$ (for $r\ge6$) are order conditions of the
composition; otherwise it only modifies the error constant (e.g.\ $\sum_s\beta_s^5\ne0$ for BM4 with $q=1$). We have verified this cancellation for non-commuting operators by a direct expansion in the free algebra generated by $\opH_1$ and $\opH_2$, up to fourth order for BM4 with $q=1$. The preservation of order relies on the time symmetry~\eqref{eq:Uq_symmetric}, not on unitarity: 
a non-symmetric approximation of the exponential, such as a truncated Taylor series, would in general reduce the order of the composition.
\end{rmk}

\begin{algorithm}[t]
\DontPrintSemicolon
\SetKwInput{Input}{Input}
\SetKwInput{Output}{Output}
\TitleOfAlgo{Combined split-operator / factored Pad\'{e} propagator}
\Input{
    State vector $\mathbf{v} = \ket{\psi(t)}$,\\
    tridiagonal operators $\opH_1, \ldots, \opH_{S-1}$ and
    diagonal operator $\opH_S = V(x,y)$,\\
    time step $\dt$,\\
    Pad\'{e} roots $\alpha_q^{(1)}, \ldots, \alpha_q^{(q)}$
    (Table~\ref{tab:Pade}),\\
    independent splitting coefficients $\beta_1, \ldots, \beta_{N_s}$
    (Table~\ref{tab:splitting})
}
\Output{$\mathbf{u} \approx U(\dt) \,\mathbf{v}$}

\BlankLine
\tcp{Full palindromic sequence: $\tilde{\beta}_s = \beta_s$ for
     $s \leq N_s$, $\tilde{\beta}_s = \beta_{2N_s+1-s}$ for $s > N_s$}
$\mathbf{u} \leftarrow \mathbf{v}$\;
\BlankLine
\For{$s = 1$ \KwTo $2N_s$}{
	\eIf{$s \leq N_s$}{
        $\tilde{\beta} \leftarrow \beta_s$\;
    }{
        $\tilde{\beta} \leftarrow \beta_{2N_s+1-s}$\;
    }
    \BlankLine
    \uIf(\tcp*[f]{$\chi(\tilde{\beta}\,\dt)$: forward order}){$s$ is odd}{
        \For{$i = 1$ \KwTo $S-1$}{
            \tcp{Factored Pad\'{e} propagator for tridiagonal $\opH_i$}
            \For{$k = 1$ \KwTo $q$}{
                $\mathbf{u} \leftarrow
                    \bigl[\Id - \alpha_q^{(k)}\rmi\,\tilde{\beta}\,\dt\,\opH_i\bigr]\,
                    \mathbf{u}$
                \tcp*{tridiagonal mv, $O(N)$}
                $\mathbf{u} \leftarrow
                    \bigl[\Id + \alpha_q^{(k)}\rmi\,\tilde{\beta}\,\dt\,\opH_i\bigr]^{-1}
                    \mathbf{u}$
                \tcp*{Thomas algorithm, $O(N)$}
            }
        }
        $\mathbf{u} \leftarrow
            \exp(-\rmi\,\tilde{\beta}\,\dt\,V)\,\mathbf{u}$
        \tcp*{diagonal: elementwise multiply, $O(N)$}
    }
    \Else(\tcp*[f]{$\chi^\star(\tilde{\beta}\,\dt)$: reversed order}){
        $\mathbf{u} \leftarrow
            \exp(-\rmi\,\tilde{\beta}\,\dt\,V)\,\mathbf{u}$
        \tcp*{diagonal: elementwise multiply, $O(N)$}
       	\For{$i = S-1, S-2, \ldots, 1$}{
            \tcp{Factored Pad\'{e} propagator for tridiagonal $\opH_i$}
            \For{$k = 1$ \KwTo $q$}{
                $\mathbf{u} \leftarrow
                    \bigl[\Id - \alpha_q^{(k)}\rmi\,\tilde{\beta}\,\dt\,\opH_i\bigr]\,
                    \mathbf{u}$
                \tcp*{tridiagonal mv, $O(N)$}
                $\mathbf{u} \leftarrow
                    \bigl[\Id + \alpha_q^{(k)}\rmi\,\tilde{\beta}\,\dt\,\opH_i\bigr]^{-1}
                    \mathbf{u}$
                \tcp*{Thomas algorithm, $O(N)$}
            }
        }
    }
}
\BlankLine
\Return{$\mathbf{u}$}
\caption{Combined propagator using the $\chi/\chi^\star$ convention
         of~\citep{McLachlan2022}. Odd stages apply
         $\chi(\tilde{\beta}\,\dt) =
         e^{-\rmi\tilde{\beta}\dt\opH_1}\cdots
         e^{-\rmi\tilde{\beta}\dt\opH_{S-1}}
         e^{-\rmi\tilde{\beta}\dt V}$;
         even stages apply $\chi^\star(\tilde{\beta}\,\dt)$ in reversed
         order. The palindromic reflection~\eqref{eq:palindrome}
         guarantees time-reversibility. Total cost:
         $O(2N_s\,q\,(S-1)\,N)$ tridiagonal operations
         $+$ $O(2N_s\,N)$ diagonal multiplications per time step.}
\label{alg:combined}
\end{algorithm}

\subsection{Numerical experiments: Two-dimensional non-separable system}
\label{sec:numerics_2D}

We consider a particle on a two-dimensional grid of $N_x\times N_y = 512\times512$ points with spacings $\dx=\delta y=0.1$~a.u., centered at the origin, with Hamiltonian $\opH = \opK_x + \opK_y + \opV$. 
The kinetic operators $\opK_x$ and $\opK_y$ are discretized by the second-order finite-difference Laplacian~\eqref{eq:H_1D} along each direction, and $\opV$ is the
diagonal matrix of the Gaussian potential
\begin{equation}
    \label{eq:V2D}
    V(x,y) = \exp\!\bigl(-\omega_x\,x^2 - \omega_y\,y^2
                          + \lambda\,x\,y\bigr)
\end{equation}
evaluated on the grid, with $\omega_x = \omega_y = 0.5$ and $\lambda = 0.3$. 
Since $V$ is not of the additive form $V_x(x)+V_y(y)$, the Hamiltonian cannot be written as a sum of commuting one-dimensional operators, and the splitting is genuinely required. 
The initial Gaussian wavepacket, centered at $(x_0, y_0) = (-5, -5)$~a.u.\ with momenta $(k_{x0}, k_{y0}) = (2, 2)$~a.u.\ and widths $(\sigma_x, \sigma_y) = (2, 2)$~a.u., is propagated to $T = 4$~a.u.

We use the decomposition~\eqref{eq:H_decomp} with $S=3$, $(\opH_1,\opH_2,\opH_3) = (\opK_x,\opK_y,\opV)$. 
With the grid points stored in lexicographic order ($x$ fastest), $\opK_x$ is tridiagonal, while $\opK_y$ couples entries separated by $N_x$.
It is therefore $N_x$-shifted and handled by Algorithm~\ref{alg:strided_thomas} with $M=N_x$ (Remark~\ref{rmk:strided}). 
Both operators thus decouple into one independent tridiagonal system per grid line, $N_y$ and $N_x$ respectively, which can be solved in parallel. 
The diagonal factor $\exp(-\rmi\,\beta_s\dt\,\opV)$ is applied exactly. 
The reference solution is again computed with \texttt{expmv} applied to the full two-dimensional Hamiltonian.

\begin{figure*}[t]
	\centering
	\includegraphics[width=1.00\textwidth]{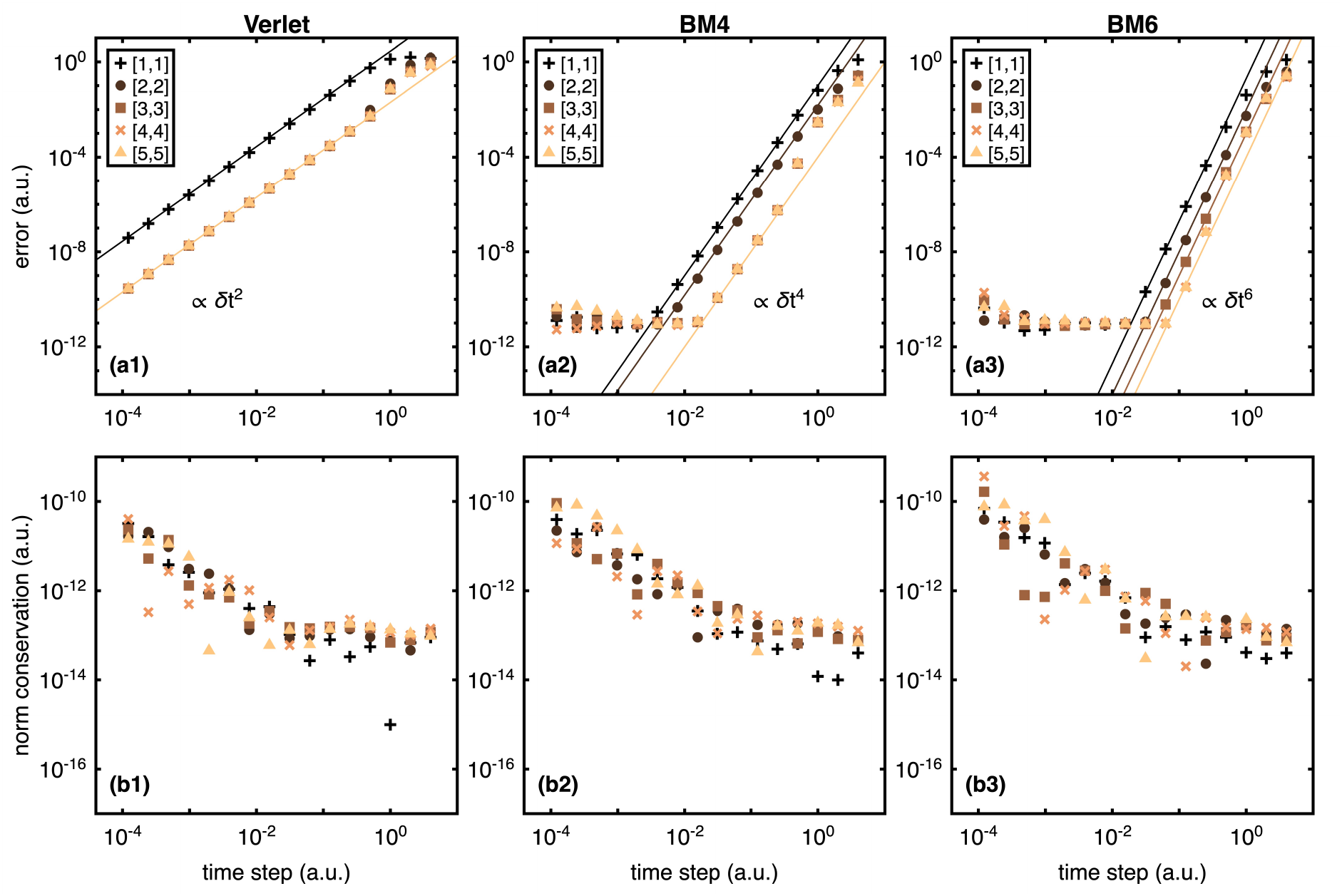}
	\caption{Two-dimensional particle with the non-separable potential~\eqref{eq:V2D}. 
	(a) $L^2$ error at time $T=4$~a.u.\ as a function of the time step $\dt$ for the diagonal Pad\'{e} families $[q,q]$, $q=1,\ldots,5$, combined with the Verlet (a1), BM4 (a2) and BM6 (a3) splitting schemes. 
	Solid lines show the scaling $\propto\dt^{r}$ of the splitting, with $r=2$ (a1), $4$ (a2) and $6$ (a3).
	Curves with $q\geq r/2+1$ coincide. 
	The error saturates at $\sim10^{-11}$, the accuracy of the \texttt{expmv} reference. 
	(b) Norm deviation at time $T$ as a function of $\dt$; its increase for small $\dt$ reflects the accumulation of rounding errors over the $N_t=T/\dt$ steps.}
	\label{fig:error_2D}
\end{figure*}

Figure~\ref{fig:error_2D}(a) shows the $L^2$ error at time $T$ as a function of $\dt$ for the Verlet, BM4 and BM6 schemes and $q=1,\ldots,5$. 
In all cases the convergence slope equals the splitting order $r$, independently of $q$, as predicted by Proposition~\ref{prop:combined}(i).
The Pad\'{e} order only affects the error constant. 
More precisely, since each time step is time-symmetric, the Pad\'{e} correction first contributes to the local error at order $\dt^{2q+1}$. 
For $q \leq r/2$ it therefore enters at the same order as the splitting error, $\dt^{r+1}$, and increases the error constant, whereas for $q \geq r/2+1$ it is of higher order and the curves coincide with that of the splitting with exact exponentials. 
This is clearly visible in Fig.~\ref{fig:error_2D}(a): the curves with $q\geq2$ coincide for Verlet, $q\geq3$ for BM4 and $q\geq4$ for BM6. 
For BM4 and BM6, the error saturates at $\sim10^{-11}$ for small $\dt$. 
As in the one-dimensional case, this floor is independent of $q$ and of the scheme, and is consistent with the accuracy of the \texttt{expmv} reference, here for a larger spectral range ($\|\opH\|\,T\approx1600$).

Figure~\ref{fig:error_2D}(b) shows the norm deviation. 
It remains at $\sim10^{-13}$ for $\dt\gtrsim10^{-2}$~a.u.\ and increases to $\sim10^{-10}$ at $\dt=10^{-4}$~a.u., for all $q$ and all schemes. 
As in Section~\ref{sec:numerics_1D}, this growth is due to the accumulation of rounding errors over the $N_t=T/\dt$ steps, the propagator being exactly unitary (Proposition~\ref{prop:combined}(ii)).

These results show that, in the presence of splitting, increasing $q$ beyond $r/2+1$ brings no further gain in accuracy, since the error is then entirely due to the splitting. 
Conversely, $q\leq r/2$ yields a larger error constant, which must be compensated by a smaller time step. 
The choice $q=r/2+1$ (i.e.\ $q=2$, $3$ and $4$ for Verlet, BM4 and BM6, respectively) therefore reaches the accuracy of the underlying splitting scheme at the lowest Pad\'{e} order.

\section{Conclusions}
\label{sec:conclusions}

These two regimes, one-dimensional propagation without splitting and multi-dimensional propagation with splitting, call for different choices of the Pad\'{e} order $q$.

In one spatial dimension (Section~\ref{sec:numerics_1D}), the Pad\'{e} approximation is the only source of temporal error. 
The convergence order is $2q$, so that increasing $q$ directly reduces the error and allows much larger time steps for a given accuracy.

In multi-dimensional problems with operator splitting (Section~\ref{sec:numerics_2D}), the convergence order is set by the splitting order $r$, independently of $q$. 
Owing to the time symmetry of the factored Pad\'{e} propagator, its contributions obey the same order conditions as the splitting itself (Proposition~\ref{prop:combined}), and the Pad\'{e} order only affects the error constant, not the scaling. 
The standard CN elementary propagator ($q=1$) therefore suffices to reach the full order of high-order splitting schemes such as BM4 and BM6, at the lowest cost per time step. 
Larger $q$ reduces the error constant, until $q \geq r/2+1$, beyond which the error coincides with that of the splitting with exact exponentials.

This also shows that the two levels of approximation are not jointly optimized: 
since the overall order is fixed by the splitting, the Pad\'{e} roots $\alpha_q^{(k)}$ represent unused degrees of freedom. 
A combined scheme in which the coefficients $\{\beta_s\}$ and $\{\alpha_q^{(k)}\}$ are tuned together could reach a target global order with fewer tridiagonal operations. 
Deriving the corresponding joint order conditions is a natural direction for future work.

\section*{Code availability}
The MATLAB implementations of the factored Pad\'{e} propagator in one and two dimensions, including the computation of the reference solutions, are available at \url{https://github.com/jonathan-dubois/High_order_Crank-Nicolson}.

\section*{Data availability}
The data underlying Figs.~\ref{fig:error_1D}--\ref{fig:error_2D} can be reproduced with the codes above.

\section*{Declaration of competing interest}
The authors declare that they have no known competing financial interests or personal relationships that could have appeared to influence the work reported in this paper.

\printcredits

\bibliographystyle{cas-model2-names}
\bibliography{cas-refs}

\end{document}